# Simulation of Permittivity and Conductivity Graded Materials for HVDC GIL for Different Voltage Forms

Hendrik Hensel, Christoph Jörgens and Markus Clemens
University of Wuppertal, Chair of Electromagnetic Theory, Rainer-Gruenter-Straße 21, 42119 Wuppertal, Germany

*Abstract*— Functionally graded materials (FGM) are applied in HVDC gas insulated lines (GIL) to control the electric field within the DC insulation system. In HVDC GIL, FGM with a spatial distribution of the electric conductivity ($\sigma$-FGM) is applied to control the electric field under DC steady state condition. However, besides DC steady state, different DC conditions occur, e.g. DC-on process, polarity reversal and lightning impulse. Under these conditions $\sigma$-FGM is not sufficient to control the electric field, since these conditions result in transient capacitive fields, where the permittivity is decisive for the electric field. In this paper, we suggest combining $\sigma$-FGM and a spatial distribution of permittivity ($\varepsilon$-FGM) in the spacer material to control the electric field around DC-GIL spacer for various DC-conditions, considering nonlinear material models for the insulating gas and the epoxy spacer. A variation of the spatial distribution of permittivity and conductivity in the spacer is investigated in this paper for an effective field reduction. The results show a reduction of the electric field intensity up to 65.8 %, when $\sigma/\varepsilon$-FGM is applied.

## I. INTRODUCTION

High voltage direct current (HVDC) has become more relevant in recent years, since HVDC is the key technology to transmit electrical energy over long distances [1]. Due to the transition to renewable electric energy, electrical energy has to be transmitted e.g. from offshore windparks to urban areas, where space efficient devices are required, since space is limited. A possible solution are gas insulated systems, which use compressed gas to gain higher breakdown voltages, to scale down the dimension of the electrical device. Gas insulated transmission lines (GIL) are used to transmit electrical energy in a space efficient way. The insulating gas used is mainly sulfur hexalfluoride ($SF_6$), which is compressed to values in the order of 1-6 bar, to obtain a dielectric strength which is up to ten times higher compared to that of air [2]. Thus, GIL require considerably less space, compared to conventional transmission devices, e.g. overhead lines. However, electric field stress distributions may occur in GIL, which can lead to partial discharge or system failure. Especially the interfaces between the insulating gas and the massive insulators in GIL, the spacers between conductor and enclosure, are sensitive to electric field stress [3]. Functionally graded materials (FGM) are used as a technique to control the electric field stress. A spatial distribution of the electric conductivity (DC) or permittivity (AC) in the spacer ($\sigma$-FGM/$\varepsilon$-FGM) have been proposed in previous works to reduce the electric field stress within the GIL [4, 5, 11, 12, 13, 14, 15, 16]. The electric conductivity is decisive for the electric field under DC steady state condition. However, transient capacitive field stress, resulting from alternating current (AC), may also appear within HVDC GIL due to e.g. polarity reversal, lightning impulse or when DC voltage is turned on/off. Under these conditions, the permittivity $\varepsilon$ of the insulating material is decisive for the electric field [5]. Hence, a combination of electric conductivity and permittivity ($\sigma/\varepsilon$-FGM) for HVDC GIL is investigated in this work, by numerical simulations. Following this introduction, the concept of $\sigma/\varepsilon$-FGM, the computation of the electric and thermal field and the nonlinear electric conductivity models of the insulating materials are presented in the next
section. The numerical simulation results for DC-on process, polarity reversal and lightning impulse with the application of $\sigma/\varepsilon$-FGM are given in section 3, followed by the conclusions in section 4.

## II. APPROACH

### A. Concept of $\sigma/\varepsilon$-FGM

The electric field within GIL consists of capacitive and resistive fields. Derived from Maxwell's equations, the electro-quasistatic complex amplitude of the potential $\varphi$ is generally determined by

$$\nabla \cdot \left( (\varepsilon - j\frac{\sigma}{2\pi f})\nabla\varphi \right) = 0, \quad (1)$$

where $\varepsilon = \varepsilon_0 \cdot \varepsilon_r$ is the permittivity, with the dielectric constant $\varepsilon_0 = 8,854 \cdot 10^{-12}$ As/Vm and the relative permittivtiy $\varepsilon_r$, $\sigma$ the electric conductivity and $f$ the frequency. Depending on the voltage condition, either the electric conductivity or the permittivity is decisive for the electric field. For DC (resistive) electric fields ($f = 0$ Hz) (1) can be converted to

$$\nabla \cdot (\sigma\nabla\varphi) = 0, \quad (2)$$

where the electric conductivity is dominant [6]. Under DC-operation condition, a slowly time varying electro-quasistatic field results in a GIL, due to space charge accumulation within the spacer and the insulating gas. The accumulation of charges causes an additional electric field, which adds to the externally applied electric field from the conductors [7]. For transient capacitive electric fields ($\varepsilon \gg \sigma/f$) the permittivity is dominant and (1) can be rewritten to

$$\nabla \cdot (\varepsilon\nabla\varphi) = 0. \quad (3)$$

However, under DC voltage capacitive and resistive fields both occur. During a DC voltage on process, the voltage level rises and results in a capacitive field, which is equivalent to the field under AC conditions. When the voltage level reaches a constant value, the slowly time varying electric fields results in a DC steady state condition, hence in a resistive field. This also occurs under

DC polarity reversal and when a lighting impulse appears under DC steady state [6]. These conditions are examined in this work, since under these conditions both $\varepsilon$-FGM and $\sigma$-FGM are required to control the electric field. The concept of FGM is to vary the permittivity and/or electric conductivity of the spacer material spatially [4, 5, 11, 12, 13, 14, 15, 16]. The areas where the spacer is under high electric fields have a higher value of permittivity/electric conductivity to reduce the electric field stress, compared to spacer parts, where high electric fields do not occur. The fabrication of FGM is realized by a controlled distribution of fillers in the spacer material, which mainly consists of epoxy resin. These fillers have a different value of permittivity/electric conductivity compared to epoxy resin to differentiate these parameters over space [3], [6].

### B. Nonlinear electric conductivity models

Under DC operation condition, the electric conductivity depends on the temperature and the electric field [6]. The heat transfer in the spacer is given by heat conduction and for the gas it can be described by convective heat transfer and heat radiation. In actual operation, the current in the conductor leads to a temperature gradient between conductor and ground. Since fixed temperature gradients are defined in the simulations, it is sufficient to solve the stationary heat conduction equation [6], [8]

$$\nabla \cdot (\lambda \nabla T) + \kappa |\vec{E}|^2 = 0. \quad (4)$$

Here, $T$ is the temperature, $\lambda$ is the thermal conductivity and $\kappa |\vec{E}|^2$ describes the heat sources in both insulation materials. Equation (1) and (4) are solved in a coupled simulation. For this, the electric conductivity needs to be determined precisely. In general, the electric conductivity of the epoxy resin material of the spacer, without considering FGM, has a nonlinear dependency on the temperature $T$ and the electric field $|\vec{E}|$, according to [4], [5], [6]. The equation to describe the electric conductivity of the spacer is given by

$$\kappa(T, |\vec{E}|) = \kappa_0 \exp\left(-\frac{W_A}{k_B T}\right) \exp(\vartheta |\vec{E}|), \quad (5)$$

where $k_B = 8.617 \cdot 10^{-5}$ eV/K is the Boltzmann constant, $W_A = 0.095$ eV the activation energy and $\kappa_0$ and $\vartheta$ are constants [5], [6]. In most works the electric conductivity of the SF$_6$ is assumed to be constant. Here, the electric conductivity is described by a model to obtain precise simulation results, which is developed in [9]. According to [9], the electric conductivity of the SF$_6$ gas depends on the electric field $|\vec{E}|$, the temperature $T$ and the gas pressure $P$. The equation which describes the electric conductivity of SF$_6$ is given by

$$\kappa(|\vec{E}|, P, T) = \kappa_{SF6} \cdot (\alpha + \beta \cdot (\gamma + |\vec{E}|/E_x)^\varsigma) \cdot (1/(\varrho + \epsilon \cdot |\vec{E}|/E_y)^\iota) \cdot \exp(\zeta \cdot P) \cdot \exp(\nu \cdot T), (6)$$

where $\kappa_{SF6}, \alpha, \beta, \gamma, E_x, \varsigma, \varrho, \epsilon, E_y, \iota, \nu$ and $\zeta$ are constants. The model allows to consider the effect different parameters have on the electric conductivity of the insulating gas, which would be neglected, if the electric conductivity was set constant.

## III. ELECTRIC FIELD SIMULATION RESULTS WITH APPLIED $\sigma/\varepsilon$-FGM

The geometry model used for the simulations is a 2D model of an axisymmetric 320 kV DC GIL with a cone-type spacer and can be depicted in Figure 1. The temperature at the ground is set at 300 K. For all simulations in this work a temperature gradient of 40 K between conductor and ground is given, to consider the influence of the temperature on the electric conductivity of the insulators and the electric field within the GIL. The gas pressure is set at 0.6 MPa for all following results. The electric conductivity models (5) and (6) are each applied for the spacer and the SF$_6$ gas. The voltage distribution within the 320 kV GIL is shown in Figure 2.

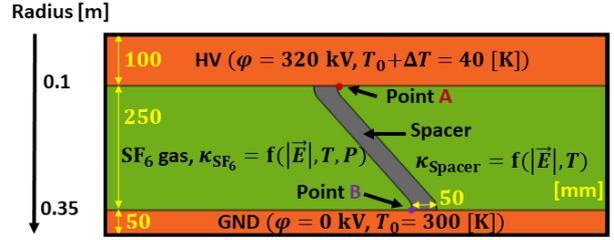

Figure 1: 2D axisymmetric geometry model of the GIL.

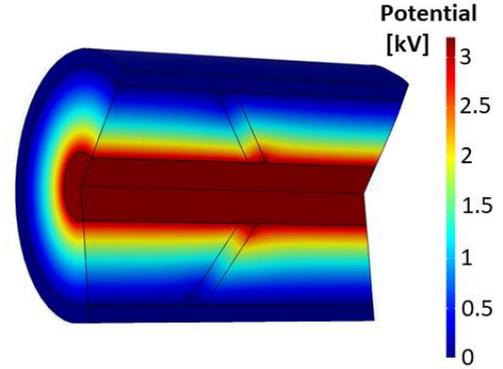

Figure 2: Potential distribution in the 320 kV GIL.

Three different DC conditions are investigated in this work: DC-on process, polarity reversal of the DC voltage and the appearance of a lightning impulse on a DC steady-state condition. Under these conditions both transient and DC fields occur and are hence convenient to investigate the benefit of $\sigma/\varepsilon$-FGM. In the following subchapters the simulation results for these conditions are presented.

### A. DC-on

The DC-on process is performed by applying the 320 kV DC voltage at the conductor within 0.01 s switching time.

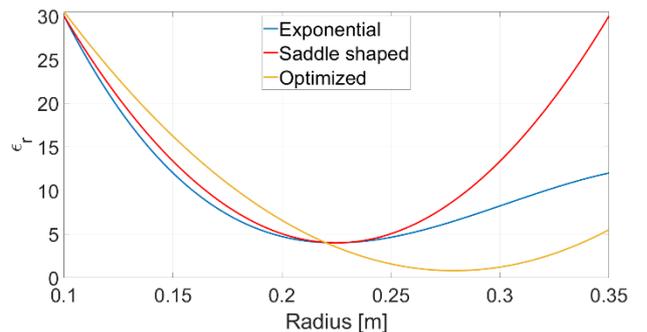

Figure 3: Three distributions of the permittivity in the spacer.

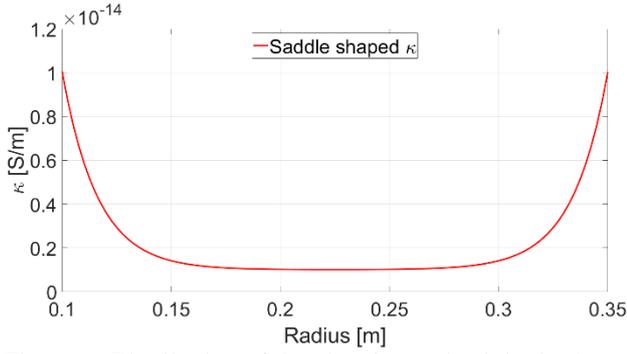

Figure 4: Distribution of the electric conductivity in the spacer.

To obtain the electric field over time, the analysis is executed by a time-domain simulation, until DC steady-state condition is reached. Three different distributions of the permittivity in the spacer material are investigated, as presented in Figure 3. These distributions are combined with a saddle shaped distribution of the electric conductivity in the spacer, depicted in Figure 4. The three different $\sigma/\varepsilon$-FGM approaches are applied on the spacer and the simulations results are compared with each other and with no application of $\sigma/\varepsilon$-FGM. Graded permittivity and electric conductivity are applied radially, which results in higher values of $\sigma/\varepsilon$ at the electrodes, where the electric field stress appears.

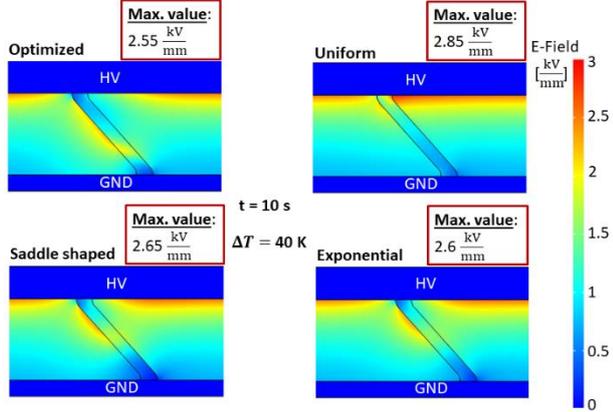

Figure 5: Distribution of the capacitive electric field under DC-on.

In Figure 5 the electric field distribution in the GIL is depicted, 10 s after the 320 kV are applied. Here, the transient capacitive field is shown with different electric field distributions, depending on the distribution of $\varepsilon$-FGM. In the case of a constant permittivity (uniform), an electric field peak occurs at the triple point between conductor, the insulating gas and the spacer. With the application of $\varepsilon$-FGM it can be seen that the electric field peak is relaxed for all different $\varepsilon$-FGM distributions. In case of $\varepsilon$-FGM, the maximum electric field is also lower, compared to the uniform case. The optimized $\varepsilon$ distribution shows the lowest maximum electric field with 2.55 kV/mm, a reduction of 10.5 % compared to the maximum electric field of 2.85 kV/mm in the uniform case.

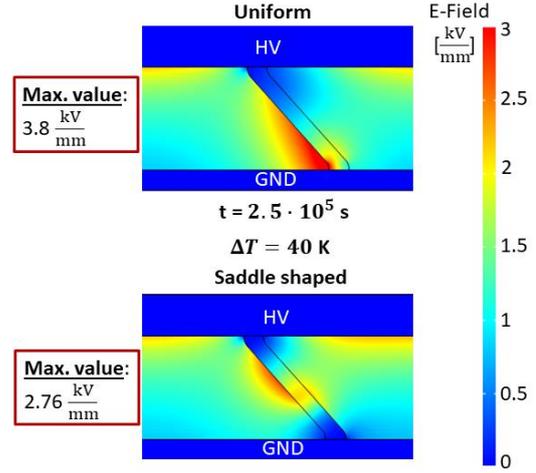

Figure 6: Distribution of the resistive electric field under DC-on.

Apart from the $\varepsilon$-FGM distributions, $\sigma$-FGM is also applied, which shows its effect after DC steady state is reached. In Figure 6 the electric field distribution is shown for t = $2.5 \cdot 10^5$ s for the uniform case and for one of the three cases where $\sigma/\varepsilon$-FGM is applied. The electric field distribution for DC steady state condition is the same in all three cases since the same electric conductivity (Figure 4) is applied. The electric field distribution in the uniform case is inversed due to the slow charge accumulation, which occurs under DC voltage. An electric field peak is seen at the triple point between the ground, the spacer and the insulating gas. With the application of $\sigma$-FGM the electric field inversion is prevented and in general is the electric field distribution relaxed. This can be also seen if the maximum electric field is examined. The maximum electric field in case of $\sigma$-FGM is 2.76 kV/mm, with a reduction of 27.3 % compared to the 3.8 kV/mm in the uniform case.

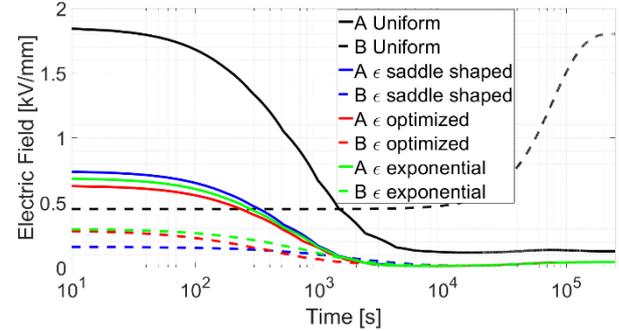

Figure 7: Electric field at the points A and B over time for different $\sigma/\varepsilon$-FGMs and a configuration without application of $\sigma/\varepsilon$-FGM (uniform).

For the illustration of the effect of a combined $\sigma/\varepsilon$-FGM, the electric field is depicted at the triple points A and B over time until DC steady state is reached, for the regarded cases in Figure 7. For all cases of application of $\sigma/\varepsilon$-FGM, during the whole process the electric field is lower at both points A (solid line) and B (dotted line) compared to the uniform spacer. Especially at point B, when DC steady state is reached, the electric field is considerably lower, since for the uniform spacer the electric field is getting inversed. Also at point A shortly after the application of the voltage the electric field is considerably higher in the uniform case. Figure 7 allows to compare the performance of the different $\sigma/\varepsilon$-FGM distributions. It can be seen that at point A the optimized $\varepsilon$-FGM shows the lowest electric field over time. At point B the optimized $\varepsilon$-FGM also indicates a very low electric field. Only for the first time period the saddle shaped $\varepsilon$-FGM shows a slightly lower electric field. It can be concluded that the optimized $\varepsilon$-FGM shows the best results. Hence, the optimized $\varepsilon$-FGM is used for further results and is compared to the uniform spacer.

### B. Polarity Reversal

When under DC voltage the polarity gets reversed, the voltage is suddenly changed, in this case from positive 320 kV to negative 320 kV in a very short switching time. Hence, both the permittivity and the electric conductivity are decisive for the electric field distribution in this process. According to (7), the electric field under polaritiy reversal $|\vec{E}_{PR}|$ can be calculated by

$$|\vec{E}_{PR}| = |2 \cdot |\vec{E}_{AC}| - |\vec{E}_{DC}||, \quad (7)$$

where $\vec{E}_{AC}$ is the capacitive electric field under AC steady state and $\vec{E}_{DC}$ is the resistive electric field under positive DC steady state [10]. The electric field distribution in the vicinity of the spacer at polarity reversal is depicted in Figure 8. The electric field is shown for a uniform spacer and for the $\sigma/\varepsilon$-FGM with the optimized permittivity distribution. In the uniform case it can be depicted that high electric field stress occurs at the triple point between spacer, gas and conductor. Compared to the $\sigma/\varepsilon$-FGM spacer, the electric field stress is highly reduced. The maximum electric field is 3.64 kV/mm and is 32.8 % reduced, compared to the maximum value of the uniform spacer, which is 5.43 kV/mm. The results indicate that the $\sigma/\varepsilon$-FGM field control technique is very useful to reduce the electric field distribution and to lower high electric field peaks, which occur under DC polarity reversal. Since both capacitive and resistive fields influence the electric field at polarity reversal, a combination of permittivity and electric conductivity FGM is necessary to control the electric field effectively.

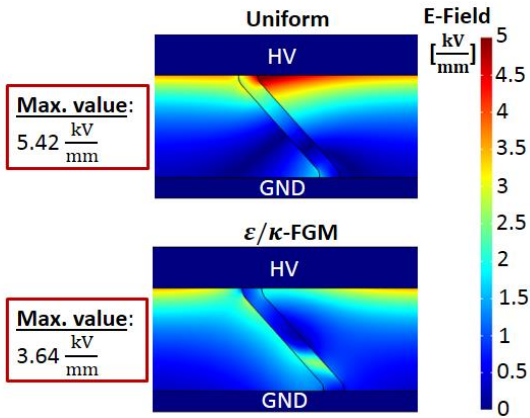

Figure 8: Electric field distribution at DC polarity reversal for a uniform spacer and a $\sigma/\varepsilon$-FGM spacer.

### C. Lightning Impulse

The 1175 kV positive lightning impulse is superimposed on 320 kV DC steady state condition. The simulation is executed through a time-domain simulation. The standard lightning impulse has a front/tail time of 1.2/50 µs, which can be seen in Figure 9, where the voltage at the conductor is illustrated over time.

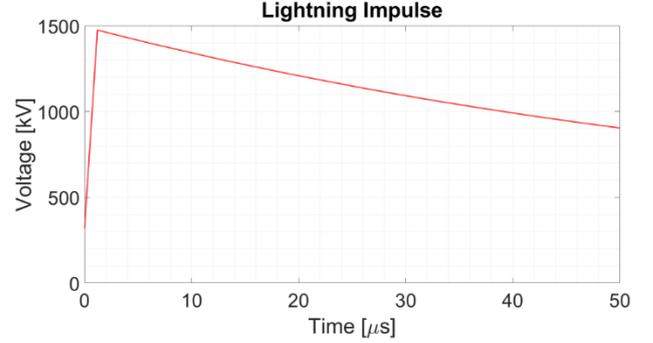

Figure 9: Superimposed 1175 kV lightning impulse on the 320 kV DC steady state.

Figure 10 shows the electric field over time after the lightning impulse was applied at point A and B for an uniform spacer and the $\sigma/\varepsilon$-FGM spacer. The electric field at point A is in case of the uniform spacer constantly higher, compared to the results for the $\sigma/\varepsilon$-FGM spacer. The electric field has its highest value right after the lightning impulse was applied, since the lightning impulse voltage decreases rapidly. The maximum value at point A for the uniform spacer is ca. 6.85 kV/mm. In terms of electric field stress, the $\sigma/\varepsilon$-FGM spacer performs clearly better. The maximum value for the $\sigma/\varepsilon$-FGM spacer is ca. 2.3 kV/mm, which is a reduction of 66.4 %, compared to the uniform case. The electric field at point B shows similar results, where the electric field is permanently lower when $\sigma/\varepsilon$-FGM is used. The maximum value of the electric field is 3.5 kV/mm and is decreased for the $\sigma/\varepsilon$-FGM spacer to 1.05 kV/mm, hence a reduction of 70 %. It can be concluded that the application of $\sigma/\varepsilon$-FGM highly reduces electric field stress under a lightning impulse.

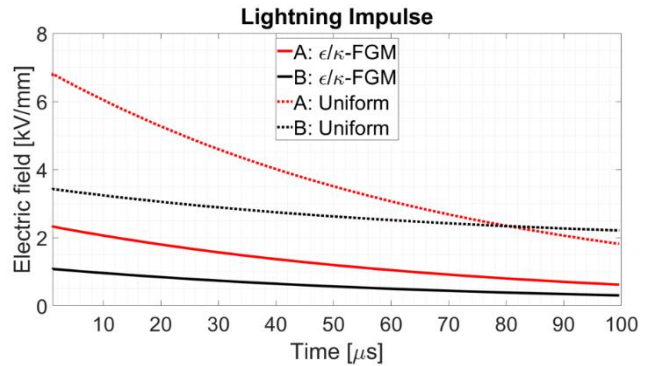

Figure 10: Electric field at point A and B for superimposed 1175 kV lightning impulse on the 320 kV DC steady state for a uniform spacer and a $\sigma/\varepsilon$-FGM spacer.

*D. Discussion of Results*

In DC-on process $\sigma/\varepsilon$-FGM proved a considerably reduced electric field intensity, especially at the triple points between gas, spacer and electrodes. The saddle shaped distribution of the electric conductivity and the optimized distribution of the permittivity in the spacer indicated the best results under the investigated distributions. The electric field stress at triple point A for the transient fields is with the application of $\sigma/\varepsilon$-FGM heavily reduced. When the electric fields develop into DC steady state, the uniform spacers shows a strongly inversed field, where the triple point B is intensely stressed. When $\sigma/\varepsilon$-FGM is applied, the inversion of the electric field is prevented, hence the electric field stress at the triple point B. Analog results are obtained for polarity reversal, where the maximum electric field within GIL is 32.8 %. When a lightning impulse overlaps with the operation DC voltage of the GIL, the GIL gets heavily stressed by the electric field, which the simulation results of the configuration without application of $\sigma/\varepsilon$-FGM show. The application of the optimized $\sigma/\varepsilon$-FGM spacer depicted excellent results, where the electric field at both sensitive triple points was considerably reduced up to 65.8 %.

## IV. Conclusion

An investigation of a combination of permittivity and electric conductivity FGM in GIL for different DC conditions was performed using transient electro-quasistatic 2D field simulation. The distribution of $\sigma/\varepsilon$-FGM in the spacer was varied to find a suitable distribution to control the electric field intensity within the GIL. To obtain precise numerical results, models for the electric conductivity for the epoxy spacer and the $SF_6$ gas were applied. Examined were DC conditions such as DC-on process, polaritiy reversal and lightning impulse.

It can be concluded that a combination of $\sigma/\varepsilon$-FGM shows very good results terms of field control techniques and could reduce the electric field significantly.

## References

[1] A. Küchler, "High voltage engineering: Fundamentals, Technology, Applications", 5[th] ed., Springer-Verlag, Germany, 2018.

[2] J. Kindesberger, C. Lederle, "Surface charge decay on Insulators in air and sulfurhexalfluoride – Part I: Simulation", IEEE Trans. Dielectr. Electr. Insul., Vol. 15, No. 4, 2008.

[3] N. Hayakawa, Y. Niyaji, H. Kojima, "Simulation on discharge inception voltage improvement of GIS spacer with permittivity graded materials ($\varepsilon$-FGM) using flexible mixture casting method", IEEE Transactions on Dielectric and Electrical Insulation, Vol. 25, No. 4, 2018, pp. 1318-1323.

[4] C. Jörgens, H. Hensel, M. Clemens, „Modeling of the electric field in high voltage direct current gas insulated transmission lines", International Conference on Dielectrics, 2022, submitted for publication.

[5] N. Hayakawa et al., "Electric field grading by functionally graded materials (FGM) for HVDC gas insulated power apparatus", IEEE Conference on Electrical Insulation and Dielectric Phenomena, Cancun, Mexico, 2018.

[6] Rachmawati et al., „Electric field simulation of permittivity and conductivity graded materials ($\varepsilon/\kappa$-FGM) for HVDC GIS spacers", IEEE Transactions on Dielectrics and Electrical Insulation, Vol. 28, No. 2, April 2022, pp. 736-744

[7] Y. Luo, et al., "Transition of the dominant charge accumulation mechanism at a Gas-solid interface under DC voltage" IET Gener. Transm. Distrib., Vol. 14, No. 15, pp. 3078-3088, 2020.

[8] H. Li, N. Zebouchi, and A. Haddad, "Theoretical and practical investigations of spacer models for future HVDC GIL/GIS applications," Proc. 21st Int. Symp. High Voltage Engineering (ISH), Vol. 2, pp. 1538–1549, 2019.

[9] H. Hensel, C. Jörgens, M. Clemens, „Numerical simulation of electric field distribution in HVDC gas insulated lines considering a novel nonlinear conductivity model for SF6", VDE Hochspannungstechnik, 2022, abstract accepted.

[10] R. Nakane et al., "Electrical insulation performance of HVDC-GIS spacer under various testing conditions", IEEE Conf. Electr. Insul. Dielect. Phenom. (CEIDP), 2017, 7-6, pp. 621-624.

[11] R. N. Hayakawa et al., „Fabrication and simulation of permittivity graded materials for electric field grading of gas insulated power apparatus", IEEE Transactions on Dielectric and Electrical Insulation, Vol. 23, No. 1, 2016, pp. 547-553.

[12] M. Kurimoto et al., "Application of functionally graded material for reducing electric field on electrode and spacer interface," IEEE Transactions on Dielectric and Electrical Insulation, Vol.17, No.1, 2010, pp. 256-263.

[13] A. Al-Gheilani, Y. Li, K. L. Wong, W. S. T. Rowe, "Electric field reduction by multi-layer functionally graded material with controlled permittivity and conductivity distribution", IEEE Conf. Electr. Insul. Dielect. Phenom., 2019, pp. 86-89.

[14] Z. Ran, "Electric field regulation of insulator interface by FGM with conductivity for superconducting-GIL", IEEE Transactions on Applied Superconductivity, Vol. 29, No. 2, 2019.

[15] B. X. Du, Z. Y. Ran, J. Li, H. C. Liang, „Novel Insulator with interfacial $\sigma$-FGM for DC compact gaseous insulated pipeline", IEEE Transactions on Dielectrics and Electrical Insulation, Vol. 26, No. 3, 2019, pp. 818-825.

[16] B. X. Du, Z. Y. Ran, J. Li, H. C. Liang, H. Yao, „Fluorinated epoxy insulator with interfacial conductivity graded material for HVDC gaseous insulated pipeline", IEEE Transactions on Dielectrics and Electrical Insulation, Vol. 27, No. 4, 2020, pp. 1305-1312.